\documentclass[runningheads, envcountsame, a4paper]{llncs}

\usepackage{epsfig}
\usepackage{graphicx}
\usepackage{amsmath}
\usepackage{amssymb}
\usepackage{mathtools}
\usepackage{algorithm}
\usepackage{enumerate}
\usepackage[noend]{algpseudocode}
\usepackage{dpr}
\usepackage{fec}
\usepackage{makecell}
\usepackage{bbm}
\usepackage{comment}
\usepackage{enumitem}
\usepackage{xcolor}
\usepackage{arydshln}
\usepackage{siunitx,tabularx,ragged2e,booktabs}
\usepackage{adjustbox}
\usepackage[misc]{ifsym}
\usepackage[hyphens]{url} 
\usepackage{microtype}

\newtheorem{assumption}{Assumption}

\DeclarePairedDelimiter\ceil{\lceil}{\rceil}

\newcommand{\rda}{\text{RDA}}
\newcommand{\proxsvrg}{\text{Prox-SVRG}}
\newcommand{\proxsg}{\text{Prox-SG}}
\newcommand{\algname}{Orthant Based Proximal Stochastic Gradient Method}
\newcommand{\algacro}{OBProx-SG{}}
\newcommand{\algacroplus}{OBProx-SG+{}}
\newcommand{\algacroall}{OBProx-SG(+){}}

\newcommand{\fashionmnist}{\text{Fashion-MNIST}}
\newcommand{\cifar}{\text{CIFAR10}}
\newcommand{\resnet}{\text{ResNet18}}

\newcommand{\mobilenet}{\text{MobileNetV1}}
\newcommand{\myeq}{\stackrel{\mathclap{\normalfont\mbox{\small def}}}{=}}

\newcommand{\prox}{\text{Prox}}
\newcommand{\objstep}{\text{Prox-SG Step}}
\newcommand{\orthantstep}{\text{Orthant Step}}

\newcommand{\ie}{\textit{i.e.}}
\newcommand{\eg}{\textit{e.g.}}
\newcolumntype{?}{!{\vrule width 1pt}}

\def \B {\mathcal{B}}

\def \R {\mathbb{R}}

\def \S {\mathbb{S}}
\def \O {\mathcal{O}}
\def \P {\mathcal{P}}

\begin{document}
	\title{Orthant Based Proximal Stochastic Gradient Method for $\ell_1$-Regularized Optimization}
	\titlerunning{OBProx-SG for $\ell_1$-Regularized Optimization}
	%
	\author{Tianyi Chen \Letter\inst{1} \and
		Tianyu Ding\inst{2} \and
		Bo Ji\inst{3} \and
		Guanyi Wang\inst{4} \and
		Yixin Shi\inst{1} \and
		Jing Tian\inst{5} \and 
		Sheng Yi\inst{1} \and
		Xiao Tu\inst{1} \and
		Zhihui Zhu\inst{6}
	}
	\authorrunning{T. Chen, T. Ding et al.}
	%
	\institute{Microsoft, Redmond, USA, \email{Tianyi.Chen@microsoft.com} \and
		Johns Hopkins University, Baltimore, USA, \email{tding1@jhu.edu} \and
		Zhejiang University, China, \email{jibo27@zju.edu.cn} \and
		Georgia Institute of Technology, Atlanta, USA, \email{gwang93@gatech.edu}\and
		University of Washington, USA, \email{jingtc20@uw.edu}\and
		University of Denver, Denver, USA, \email{zhihui.zhu@du.edu}. 
	}

	\maketitle              
	
	\setcounter{footnote}{0}

	\begin{abstract}
		Sparsity-inducing regularization problems are ubiquitous in machine learning applications, ranging from feature selection to model compression.  In this paper, we present a novel stochastic method -- Orthant Based Proximal Stochastic Gradient Method (\algacro{}) -- to solve perhaps the most popular instance, \ie, the $\ell_1$-regularized problem. The \algacro{} method contains two steps: (i) a proximal stochastic gradient step  to predict a support cover of the solution; and (ii) an orthant step to aggressively enhance the sparsity level via orthant face projection. Compared to the state-of-the-art methods, \eg, \proxsg, \rda{} and \proxsvrg, the \algacro{} not only converges comparably in both convex and non-convex scenarios, but also promotes the sparsity of the solutions substantially. Particularly, on a large number of convex problems, OBProx-SG  outperforms the existing methods comprehensively in the aspect of sparsity exploration and objective values. Moreover, the experiments on non-convex deep neural networks, \eg, MobileNetV1 and ResNet18, further demonstrate its superiority by generating the solutions of much higher sparsity without sacrificing generalization accuracy, which further implies that~\algacro{} may achieve significant memory and energy savings. The source code is available at {\color{magenta}\url{https://github.com/tianyic/obproxsg}}.
		
		\keywords{Stochastic Learning \and Sparsity \and Orthant Prediction.}
	\end{abstract}

	\section{Introduction}\label{sec.introduction}
	
	Plentiful tasks in machine learning and deep learning require formulating and solving particular optimization problems~\cite{bradley1977applied,dixit1990optimization}, of which the solutions may not be unique. From the perspective of the application, people are usually interested in a subset of the solutions with certain properties. A common practice to address the issue is to augment the objective function by adding a regularization term~\cite{tikhonov1977}. 	
	One of the best known examples is the sparsity-inducing regularization, which encourages highly sparse solutions (including many zero elements). Besides, such regularization typically has shrinkage effects to reduce the magnitude of the solutions~\cite{tibshirani1996regression}. Among the various ways of introducing sparsity, the $\ell_1$-regularization is perhaps the most popular choice. Its utility has been demonstrated ranging from improving the interpretation and accuracy of model estimation~\cite{riezler2004incremental,sra2011fast} to compressing heavy model for efficient inference~\cite{cheng2017survey,han2015deep}. 
	
	In this paper, we propose and analyze a novel efficient stochastic method to solve the following large-scale $\ell_1$-regularization problem
	\begin{equation}\label{prob.x}
	\minimize{x\in \R^n}\ \Big\{F(x)\ \myeq\ \underbrace{\frac{1}{N}\sum_{i=1}^N f_i(x)}_{f(x)}+\lambda\norm{x}_1\Big\},
	\end{equation}
	where $\lambda>0$ is a weighting term to control the level of sparsity in the solutions, and $f(x)$ is the raw objective function. We pay special interests to the $f(x)$ as the average of numerous $N$ continuously differentiable instance functions $f_i : \R^n\rightarrow\R$, such as the loss functions  measuring the deviation from the observations in various data fitting problems. A larger $\lambda$ typically results in a higher sparsity while sacrifices more on the bias of model estimation, hence $\lambda$ needs to be carefully fine-tuned to achieve both low $f(x)$ and high sparse solutions. Above formulation is widely appeared in many contexts, including convex optimization, \eg, LASSO, logistic regression and elastic-net formulations~\cite{tibshirani1996regression,zou2005regularization}, and non-convex problems such as deep neural networks~\cite{zaremba2014recurrent,zeiler2013stochastic}.
	
	Problem~\eqref{prob.x} has been well studied in deterministic optimization with various methods that capable of returning solutions with both low objective value and high sparsity under proper $\lambda$. Proximal methods are classical approaches to solve the structured non-smooth optimization problems with the formulation~\eqref{prob.x}, including the popular proximal gradient method (Prox-FG) and its variants, \eg,  ISTA and FISTA~\cite{beck2009fast}, in which only the first-order derivative information is used. They have been proved to be quite useful in practice because of their simplicity. Meanwhile, first-order methods are limited due to the local convergence rate and lack of robustness on ill-conditioned problems, which can often be overcome by employing the second-order derivative information as is used in proximal-Newton methods~\cite{lee2012proximal,yuan2012improved}. However, when $N$ is enormous, a straightforward computation of the full gradients or Hessians could be prohibitive because of the costly evaluations over all $N$ instances. Thus, in modern large-scale machine learning applications, it is inevitable to use stochastic methods that operate on a small subset of above summation to economize the computational cost at every iteration. 
	
	Nevertheless, in stochastic optimization, the studies of $\ell_1$-regularization~\eqref{prob.x} become somewhat limited. In particular, the existing state-of-the-art stochastic algorithms rarely achieve both fast convergence and highly sparse solutions simultaneously due to the stochastic nature~\cite{xiao2010dual}. Proximal stochastic gradient method (Prox-SG)~\cite{duchi2009efficient} is a natural extension of Prox-FG by using a mini-batch to estimate the full gradient. However, there are two potential drawbacks of Prox-SG: (i) the lack of exploitation on the certain problem structure, \eg, the $\ell_1$ regularization \eqref{prob.x}; (ii) the slower convergence rate than Prox-FG due to the variance introduced by random sampling. To exploit the regularization structure more effectively (produce sparser solutions), regularized dual-averaging method (\rda)~\cite{xiao2010dual} is proposed by extending the simple dual averaging scheme in~\cite{nesterov2009primal}.
	The key advantages of \rda{} are to utilize the averaged accumulated gradients of $f(x)$ and an aggressive coefficient of the proximal function to achieve a more aggressive truncation mechanism than~\proxsg{}. As a result, in convex setting, \rda{} usually generates much sparser solutions than that by \proxsg{} in solving \eqref{prob.x}  but typically has slower convergence. On the other hand, to reduce the variance brought by the stochastic approximation, proximal stochastic variance-reduced gradient method (\proxsvrg{})~\cite{xiao2014proximal} is developed based on the well-known variance reduction technique SVRG developed in~\cite{johnson2013accelerating}. 	\proxsvrg{} has both capabilities of decent convergence rate and sparsity exploitation in convex setting, while its per iteration cost is much higher than other approaches due to the calculation of full gradient for achieving the variance reduction. 
	
	The above mentioned \proxsg{}, \rda{} and~\proxsvrg{} are valuable state-of-the-art stochastic algorithms with apparent strength and weakness. \rda{} and~\proxsvrg{} are derived from proximal gradient methods, and make use of different averaging techniques cross all instances to effectively exploit the problem structure. Although they explore sparsity well in convex setting, the mechanisms may not perform as well as desired in non-convex formulations~\cite{defazio2019ineffectiveness}. Moreover, observing that the proximal mapping operator is applicable for any non-smooth penalty function, this generic operator may not be sufficiently insightful if the regularizer satisfies extra properties. In particular, the non-smooth $\ell_1$-regularized problems of the form~\eqref{prob.x} degenerate to a smooth reduced space problem if zero elements in the solution are correctly identified. 
	
	This observation has motivated the exploitation of orthant based methods, a class of deterministic second-order methods that utilizes the particular structure within the $\ell_1$-regularized problem~\eqref{prob.x}. During the optimization, they predict a sequence of orthant faces, and minimize smooth quadratic approximations to~\eqref{prob.x} on those orthant faces until a solution is found~\cite{andrew2007scalable,chen2017reduced,keskar2015second,chen2018fast}. Such a process normally equips with second-order techniques to yield superlinear convergence towards the optimum, and introduces sparsity by Euclidean projection onto the constructed orthant faces. Orthant based methods have been demonstrated competitiveness in deterministic optimization to proximal methods~\cite{chen2017reduced,chen2018farsa,keskar2015second}. In contrast, related prior work in stochastic settings is very rare, perhaps caused by the expensive and non-reliable orthant face selection under randomness. 
	
	\paragraph{\textbf{Our Contributions.}} In this paper, we propose an~\algname{} (\algacro) by capitalizing on the advantages of orthant based methods and Prox-SG, while avoiding their disadvantages. Our~\algacro{} is efficient, promotes sparsity more productively than others, and converges well in both practice and theory. Specifically, we have the following contributions.
	
	\begin{itemize}[leftmargin=*]
		\item We provide a novel stochastic algorithmic framework that utilizes Prox-SG Step and reduced space Orthant Step to effectively solve problem~\eqref{prob.x}. Comparing with the existing stochastic algorithms, it exploits the sparsity significantly better by combining the moderate truncation mechanism of \proxsg{} and an aggressive orthant face projection under the control of a switching mechanism. The switching mechanism is specifically established in the stochastic setting, which is simple but efficient, and performs quite well in practice. Moreover, we present the convergence characteristics under both convex and non-convex formulations, and provide analytic and empirical results to suggest the strategy of the inherent switching hyperparameter selection.
		
		\item We carefully design the Orthant Step for stochastic optimization in the following aspects: (i) it utilizes the sign of the previous iterate to select an orthant face, which is more efficient compared with other strategies that involve computations of (sub)-gradient in the deterministic orthant based algorithms~\cite{andrew2007scalable,keskar2015second}; (ii) instead of optimizing with second-order methods, only the first-order derivative information is used to exploit on the constructed orthant face. 
		
		\item Experiments on both convex (logistic regression) and non-convex (deep neural networks) problems show that OBProx-SG usually outperforms  the other state-of-the-art methods comprehensively in terms of the sparsity of the solution, final objective value, and runtime. Particularly, in the popular deep learning applications, without sacrificing generalization performance, the solutions computed by~\algacro{} usually possess multiple-times higher sparsity than those searched by the competitors.
		
	\end{itemize}
	
	\section{The \algacro\ Method}\label{sec.algorithm}
	
	To begin, we summarize the proposed~\algname{} (\algacro) in Algorithm~\ref{alg:main.x.outline}. In a very high level, it proceeds one of the two subroutines at each time, so called \objstep{} (Algorithm~\ref{alg:main.x.prox_sg_step}) and \orthantstep{}  (Algorithm~\ref{alg:main.x.orthantstep}). There exist two switching parameters $N_\mathcal{P}$ and $N_\mathcal{O}$ that control how long we are sticking to each step and when to switch to the other.  We will see that the switching mechanism (choices of $N_\mathcal{P}$ and $N_\mathcal{O}$) is closely related to the convergence of OBProx-SG and the sparsity promotions. But we defer the detailed discussion till the end of this section, while first focus our attention on the \objstep{} and \orthantstep{}.
	
	\begin{algorithm}[h!]
		\caption{Outline of \algacro{} for solving \eqref{prob.x}.}
		\label{alg:main.x.outline}
		\begin{algorithmic}[1]
			\State \textbf{Input:} $x_0\in\mathbb{R}^n$, $ \alpha_0\in(0,1)$, and $ \{N_\mathcal{P}, N_\mathcal{O}\}\subset\mathbb{Z}^+ $.
			\For{$k = 0,1,2,\dots$}
			\State \textbf{Switch}  Prox-SG Step or Orthant Step by Algorithm~\ref{alg:main.x.switch}. 
			\If{\objstep{} is selected} \label{line:switch_prox_sg_step}
			\State 
			Compute the Prox-SG Step update: 
			
			\ \ \ \ $x_{k+1}\leftarrow \text{Prox-SG}(x_k,\alpha_k)$ by Algorithm~\ref{alg:main.x.prox_sg_step}.
			
			\ElsIf{Orthant Step is selected}  \label{line:switch_orthant_step}
			
			\State Compute the Orthant Step update:
			
			\ \ \ \ $x_{k+1}\leftarrow\text{Orthant}(x_k,\alpha_k)$ 
			by Algorithm \ref{alg:main.x.orthantstep}.
			\EndIf
			\State Update $\alpha_{k+1}$ given $\alpha_k$ according to some rule.
			\EndFor
		\end{algorithmic}
	\end{algorithm}
	
	\begin{algorithm}[t]
		\caption{Prox-SG Step.}
		\label{alg:main.x.prox_sg_step}
		\begin{algorithmic}[1]
			\State \textbf{Input:} Current iterate $x_k$, and step size $ \alpha_k$. 
			\State Compute the stochastic gradient of $f$ on $ \mathcal{B}_k $
			\begin{equation}
			\nabla f_{\mathcal{B}_k}(x_k)\leftarrow\frac{1}{|\mathcal{B}_k|}\sum_{i\in \mathcal{B}_k}\Grad f_i(x_k).
			\end{equation}\label{line:g_t_estimate_prox_sg}
			\State Compute
			$x_{k+1}\leftarrow\prox_{\alpha_k\lambda\|\cdot\|_1}\left(x_k-\alpha_k\nabla f_{\mathcal{B}_k}(x_k)\right)$ \label{line:prox}.
			\State \textbf{Return } $x_{k+1}$.
		\end{algorithmic}
	\end{algorithm}

	\paragraph{\textbf{Prox-SG Step.}}
	In Prox-SG step, the algorithm performs one iteration of standard proximal stochastic gradient step to approach a solution of~\eqref{prob.x}. Particularly, at $k$-th iteration, we sample a mini-batch $\mathcal{B}_k$ to make an unbiased estimate of the full gradient of $f$ (line~\ref{line:g_t_estimate_prox_sg}, Algorithm~\ref{alg:main.x.prox_sg_step}). Then we utilize the following proximal mapping to yield next iterate as 
	\begin{equation}\label{eq:proxmapping}
	\begin{split}
	x_{k+1}&=\prox_{\alpha_k\lambda\norm{\cdot}_1}(x_k-\alpha_k \Grad f_{\mathcal{B}_k}(x_k))\\
	&=\argmin_{x\in \mathbb{R}^{n}}\ \frac{1}{2\alpha_k}\norm{x-(x_k-\alpha_k \Grad f_{\mathcal{B}_k}(x_k))}_2^2+ \lambda\norm{x}_1.
	\end{split}
	\end{equation}
	\noindent
	It is known that the above subproblem~\eqref{eq:proxmapping} has a closed form solution \cite{beck2009fast}. Denote the trial iterate $\widehat{x}_{k+1}:=x_k-\alpha_k \Grad f_{\mathcal{B}_k}(x_k)$, then $x_{k+1}$ is computed efficiently as:
	\begin{equation}\label{eq:shrink}
	\begin{split}
	[x_{k+1}]_i&=\begin{cases}
	[\widehat{x}_{k+1}]_i-\alpha_k\lambda, & \text{if}\  [\widehat{x}_{k+1}]_i > \alpha_k\lambda\\
	[\widehat{x}_{k+1}]_i+\alpha_k\lambda, & \text{if}\  [\widehat{x}_{k+1}]_i < - \alpha_k\lambda\\
	0, &  \text{otherwise}
	\end{cases}.
	\end{split}
	\end{equation}
	In~\algacro{}, \objstep{} generally serves as a globalization mechanism to guarantee convergence and predict a cover of supports (non-zero entries) in the solution. But it alone is insufficient to exploit the sparsity structure because of the relatively moderate truncation mechanism effected in a small projection region, i.e., the trial iterate $\widehat{x}_{k+1}$ is projected to zero only if it falls into $[-\alpha_k\lambda, \alpha_k\lambda]$. Our remedy here is to combine it with our Orthant Step, which exhibits an aggressive sparsity promotion mechanism while still remains efficient.
	
	\begin{algorithm}[t]
		\caption{\orthantstep{}.}
		\label{alg:main.x.orthantstep}
		\begin{algorithmic}[1]
			\State \textbf{Input:} Current iterate $x_k$, and step size $ \alpha_k$.
			\State Compute the stochastic gradient of {${\widetilde{F}}$} on $ \mathcal{B}_k $
			\begin{equation}\label{line:g_t_estimate_orthant}
			\Grad \widetilde{F}_{\mathcal{B}_k}(x_k)\gets \frac{1}{|{\mathcal{B}_k}|}\sum_{i\in \mathcal{B}_k}\Grad \widetilde{F}_i(x_k)
			\end{equation}
			\State Compute $x_{k+1}\leftarrow\proj_{\mathcal{O}_k}(x_k-\alpha_k\Grad \widetilde{F}_{\mathcal{B}_k}(x_k))$.\label{line:y_tp1_euclidean_proj}
			\State \textbf{Return } $x_{k+1}$.
		\end{algorithmic}
	\end{algorithm}
	
	\paragraph{\textbf{Orthant Step.}}  Since the fundamental to Orthant Step is the manner in which we handle the zero and non-zero elements, we define the following index sets for any $x\in\mathbb{R}^n $:
	\begin{equation}\label{def:I_set}
	\mathcal{I}^0(x) := \{i: [x]_i=0\},\ \mathcal{I}^+(x) :=\{i: [x]_i>0\},\ \mathcal{I}^-(x) :=\{i: [x]_i<0\}.
	\end{equation}
	Furthermore, we denote the non-zero indices of $x$ by $\mathcal{I}^{\neq0}(x):=\mathcal{I}^+(x)\cup\mathcal{I}^-(x)$. To proceed, we define the orthant face $\mathcal{O}_k$ that $x_k$ lies in to be 
	\begin{equation}\label{def:orthant_face}
	\mathcal{O}_k:=\{x\in\mathbb{R}^n : \sign([x]_i) = \sign([x_k]_i) \text{ or } [x]_i=0, 1\le i\le n\}
	\end{equation}
	such that $x\in\mathcal{O}_k$ satisfies: (i) $[x]_{\mathcal{I}^0(x_k)}=0$; (ii) for $i\in\mathcal{I}^{\neq0}(x_k)$,  $[x]_i$ is either 0 or has the same sign as $[x_k]_i$.
	
	The key assumption for Orthant Step is that an optimal solution $x^*$ of problem \eqref{prob.x} inhabits $\mathcal{O}_k$, \ie, $x^*\in\mathcal{O}_k$. In other words, the orthant face $\mathcal{O}_k$ already covers the support (non-zero entries) of $x^*$.
	Our goal becomes now minimizing $F(x)$ over $\mathcal{O}_k$,
	\ie, solving the following subproblem:
	\begin{equation}\label{prob.orthant_sub_problem}
	x_{k+1}=\argmin_{x\in\mathcal{O}_k}\ F(x)=f(x)+\lambda\norm{x}_1.
	\end{equation}
	By the definition of $\mathcal{O}_k$, we know $[x]_{\mathcal{I}^0(x_k)}\equiv 0$ are fixed, and only the entries of $[x]_{\mathcal{I}^{\neq0}(x_k)}$ are free to move.  Hence, ~\eqref{prob.orthant_sub_problem} is essentially a reduced space optimization problem. Observing that for any $x\in\mathcal{O}_k$, $F(x)$ can be written precisely as a smooth function $\widetilde{F}(x)$ in the form 
	\begin{equation}
	F(x)\equiv \widetilde{F}(x)\coloneqq f(x)+\lambda \sign(x_k)^Tx,
	\end{equation}
	therefore
	~\eqref{prob.orthant_sub_problem} is equivalent to the following smooth problem
	\begin{equation}\label{prob.orthant_equiv_sub_problem}
	x_{k+1}=\argmin_{x\in\mathcal{O}_k}\ \widetilde{F}(x).
	\end{equation}
	
\begin{figure}[]
		\centering
		\includegraphics[width=0.5\textwidth]{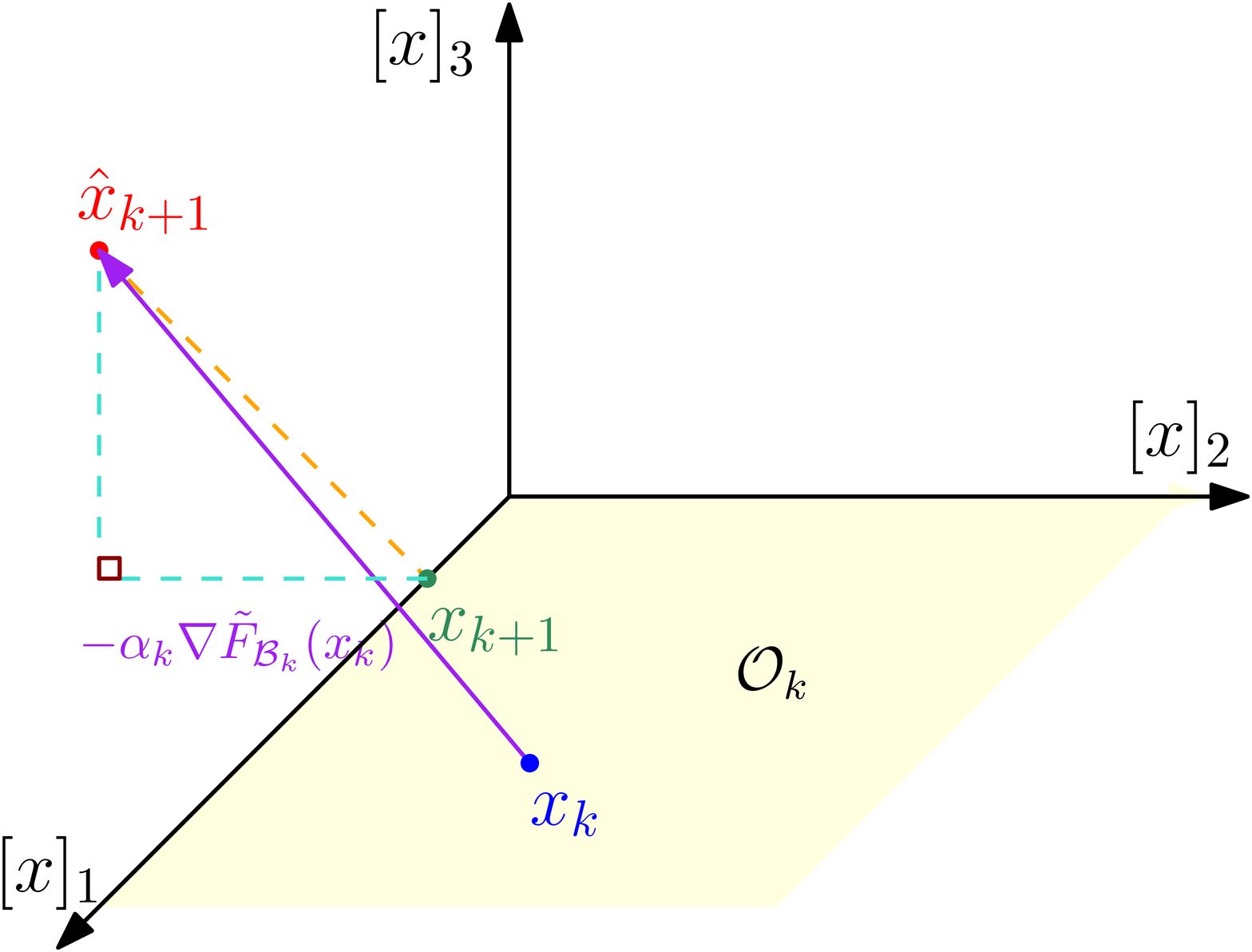}
		\includegraphics[width=0.49\textwidth]{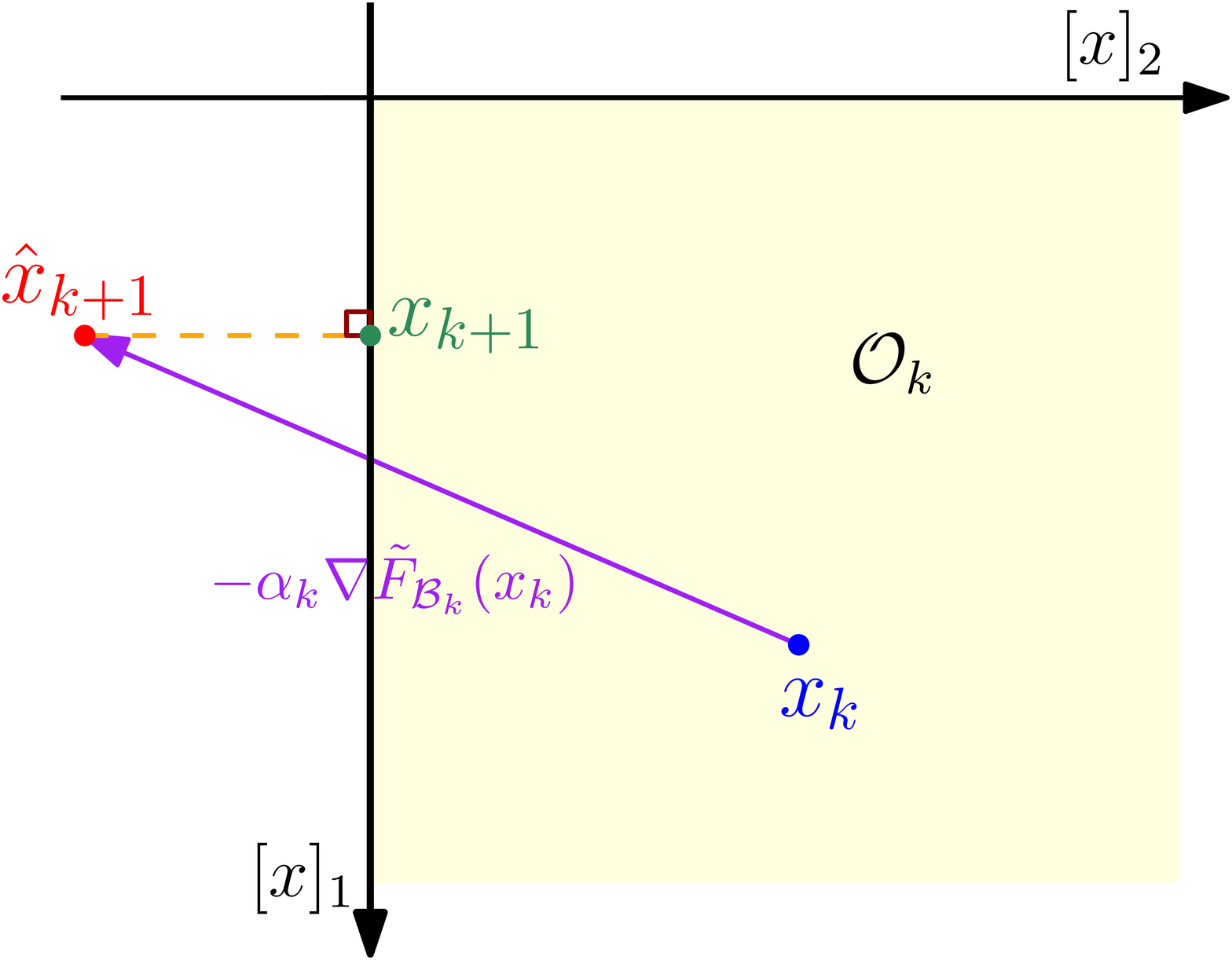}
		\caption{Illustration of Orthant Step with projection in \eqref{def:proj}, where $\mathcal{O}_k = \{x\in \mathbb{R}^3: [x]_1\ge0, [x]_2 \ge 0, [x]_3 = 0\}$. (L): 3D view. (R): top view.}
		\label{figure:proj_euclidean}
	\end{figure}

\begin{figure}[]
		\centering
		\includegraphics[width=0.7\textwidth]{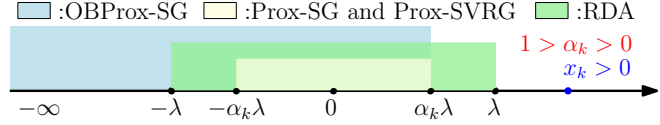}
		\caption{Projection regions of different methods for 1D case at $x_k>0$.}
		\label{figure:proj_regions}
	\end{figure}

	A direct way for solving problem~\eqref{prob.orthant_equiv_sub_problem} is the  projected stochastic gradient descent method, as stated in Algorithm~\ref{alg:main.x.orthantstep}. It performs one iteration of stochastic gradient descent (SGD) step combined with projections onto the orthant face $\mathcal{O}_k$. At $k$-th iteration, a mini-batch $\mathcal{B}_k$ is sampled, and is used to approximate the full gradient $\nabla\widetilde{F}(x_k)$ by the unbiased estimator $\Grad \widetilde{F}_{\mathcal{B}_k}(x_k)$ (line~\ref{line:g_t_estimate_prox_sg} , Algorithm~\ref{alg:main.x.orthantstep}). The standard SGD update computes a trial point $\widehat{x}_{k+1}= x_k-\alpha_k \Grad \widetilde{F}_{\mathcal{B}_k}(x_k)$, which is then passed into a projection operator $\proj_{\mathcal{O}_k}(\cdot)$ defined as
	\begin{equation}\label{def:proj}
	\left[\proj_{\mathcal{O}_k}(z)\right]_i\coloneqq\bigg\{
	\begin{array}{ll}
	[z]_i &\quad \text{if}\ \sign([z]_i) = \sign([x_k ]_i)\\ 
	0 &\quad \text{otherwise}
	\end{array}.
	\end{equation}
	Notice that $\proj_{\mathcal{O}_k}(\cdot)$ is an Euclidean projector, and ensures that the trial point $\widehat{x}_{k+1}$ is projected back to the current orthant face $\mathcal{O}_k$ if it happens to be outside, as illustrated in Figure~\ref{figure:proj_euclidean}. In the demonstrated example, the next iterate $x_{k+1}=\proj_{\mathcal{O}_k}(\widehat{x}_{k+1})$ turns out to be not only a better approximated solution but also sparser compared with $x_k$ since $[x_{k+1}]_2=0$ after the projection, which suggests the power of Orthant Step in sparsity promotion. In fact, compared with~\proxsg{}, the orthant-face projection~\eqref{def:proj} is a more aggressive sparsity truncation mechanism. Particularly, \orthantstep{} enjoys a much larger projection region to map a trial iterate to zero comparing with other stochastic algorithms. Consider the 1D example in Figure~\ref{figure:proj_regions}, where $x_k > 0$, it is clear that the projection region of \orthantstep{} $(-\infty, \alpha_k\lambda]$ is a superset of that of Prox-SG and Prox-SVRG $[-\alpha_k\lambda, \alpha_k\lambda]$, and it is apparently larger than that of RDA. 

	In practice, by taking advantage of the fact that \eqref{prob.orthant_equiv_sub_problem} is a reduced space problem, \ie, $[x_{k+1}]_{\mathcal{I}^0(x_k)}\equiv0$, we only need to store a small part of stochastic gradient information $[\nabla\widetilde{F}(x_k)]_{\mathcal{I}^{\neq0}(x_k)}$, and compute the projection of  $[\widetilde{x}_{k+1}]_{\mathcal{I}^{\neq0}(x_k)}$. This makes the whole procedure even more efficient when $\big|\mathcal{I}^{\neq0}(x_k)\big|\ll n$.
	
	We emphasize that Orthant Step is one of the keys to the success of our proposed OBProx-SG method in terms of sparsity exploration. It is originated from the orthant based methods in deterministic optimization, which normally utilize second-order information. When borrowing the idea, we make the selection of orthant face $\mathcal{O}_k$ more effective by looking at the sign of the previous iterate (see \eqref{def:orthant_face}). Then, we make use of a stochastic version of the projected gradient method in solving the subproblem \eqref{prob.orthant_equiv_sub_problem} to introduce sparsity aggressively. As a result, \orthantstep{} always makes rapid progress to the optimality, while at the same time promotes sparse solutions dedicatedly.
	

	\paragraph{\textbf{Switching Mechanism.}} 
	To complete the discussion of the~\algacro{} framework, we now explain how we select \proxsg{} or~\orthantstep{} at each iteration, which is crucial in generating accurate solutions with high sparsity. A popular switching mechanism for deterministic multi-routine  optimization algorithms utilizes the optimality metric of each routine (typically the norm of (sub)gradient)~\cite{chen2017reduced,chen2018farsa}. However, in stochastic learning, this approach does not work well in general due to the additional computation cost of optimality metric and the randomness that may deteriorate the progress of sparsity exploration as numerically illustrated in Appendix C.
	
	
	\begin{algorithm}[]
		\caption{Switching Mechanism.}
		\label{alg:main.x.switch}
		\begin{algorithmic}[1]
			\State \textbf{Input:} $k, N_\mathcal{P}, N_\mathcal{O}$.
			\If {$\text{mod} (k,N_\mathcal{P}+N_\mathcal{O})< N_\mathcal{P}$}
			\State \textbf{Return} \objstep{}  is selected. 
			\Else
			\State \textbf{Return} \orthantstep{} is selected. 
			\EndIf
		\end{algorithmic}
	\end{algorithm}
	
	To address this issue, we specifically establish a simple but efficient switching mechanism consisting of two hyperparameters $N_\mathcal{P}$ and $N_\mathcal{O}$, which performs quite well in both practice and theory. As stated in Algorithm~\ref{alg:main.x.switch}, $N_\mathcal{P}$ ($N_\mathcal{O}$) controls how many consecutive iterations we would like to spend doing \objstep{} (\orthantstep), and then switch to the other. \algacro{} is highly flexible to different choices of $N_\mathcal{P}$ and $N_\mathcal{O}$. For example, an alternating strategy between one~\objstep{} and one~\orthantstep{} corresponds to set $N_\mathcal{P}=N_\mathcal{O}=1$,~and a strategy of first doing several Prox-SG Steps then followed by Orthant Step all the time corresponds to set $N_\mathcal{P}<\infty, N_\mathcal{O}=\infty$.~A larger $N_\mathcal{P}$ helps to predict a better orthant face $\mathcal{O}_k$ which hopefully covers the support of $x^*$, and a larger $N_\mathcal{O}$ helps to explore more sparsity within $\mathcal{O}_k$.


	As we will see in Section~\ref{sec:convergence}, convergence of Algorithm~\ref{alg:main.x.outline} requires either doing \objstep{}  infinitely many times ($N_{\mathcal{P}}\le\infty, N_{\mathcal{O}}<\infty$), or doing finitely many Prox-SG Steps followed by infinitely many Orthant Steps ($N_{\mathcal{P}}<\infty, N_{\mathcal{O}}=\infty$) given the support of $x^*$ has been covered by some $\mathcal{O}_k$. In practice, without knowing $x^*$ ahead of time, we can always start by employing~\objstep{} $N_{\mathcal{P}}$ iterations, followed by running~\orthantstep{} $N_{\mathcal{O}}$ iterations, then repeat until convergence. Meanwhile, experiments in Section~\ref{sec:exp} show that first performing \objstep{} sufficiently many times then followed by running~\orthantstep{} all the time usually produces even slightly better solutions. Moreover, for the latter case, 
	a bound for $N_\mathcal{P}$ is provided in Section~\ref{sec:convergence}. For simplicity, we refer the~\algacro{} under ($N_{\mathcal{P}}<\infty, N_{\mathcal{O}}=\infty$) as~\algacroplus{} throughout the remainder of this paper.

	
	We end this section by giving empirical suggestions of setting $N_{\mathcal{P}}$ and $N_\mathcal{O}$. Overall, in order to obtain accurate solutions of high sparsity, we highly recommend to start OBProx-SG with \proxsg{} Step and ends with Orthant Step. Practically, employing finitely many Prox-SG Steps followed by sticking on Orthant Steps ($N_{\mathcal{P}}<\infty, N_{\mathcal{O}}=\infty$) until the termination, is more preferable because of its attractive property regarding maintaining the progress of sparsity exploration. In this case, although the theoretical upper bound of $N_\mathcal{P}$ is difficult to be measured, we suggest to keep running~\proxsg{} Step until reaching some acceptable evaluation metrics e.g., objectives or validation accuracy, then switch~\orthantstep{} to promote sparsity. 
	
	\section{Convergence Analysis} \label{sec:convergence}
	
	In this section,we give a convergence analysis of our proposed~\algacro{} and OBProx-SG+, referred as \algacroall{} for simplicity. Towards that end, we first make the following assumption.

	
	\begin{assumption}\label{assumption}
		The function $f:\mathbb{R}^n\to \mathbb{R}$ is continuously differentiable, and bounded below on the compact level set $\mathcal{L}\coloneqq\{x\in\mathbb{R}^n: F(x)\leq F(x_0)\}$, where $ x_0 $ is the initialization of Algorithm~\ref{alg:main.x.outline}. The stochastic gradient $\Grad f_{\mathcal{B}_k}$ and $\Grad \widetilde{F}_{\mathcal{B}_k}$  evaluated on the mini-batch $ \mathcal{B}_k $ are Lipschitz continuous on the level set $\mathcal{L}$ with a shared Lipschitz constant $L$ for all $\mathcal{B}_k$. The gradient $\Grad \widetilde{F}_{\mathcal{B}_k}(x)$ is uniformly bounded over $\mathcal{L}$, i.e., there exists a $M<\infty$ such that $\norm{\Grad \widetilde{F}_{\mathcal{B}_k}(x)}_2\leq M$.
	\end{assumption}
	
	Remark that many terms in Assumption~\ref{assumption} appear in numerical optimization literatures~\cite{chen2017reduced,xiao2014proximal,yang2019stochastic}. Let $x^*$ be an optimal solution of problem \eqref{prob.x}, $F^*$ be the minimum, and $\{x_k\}_{k=0}^\infty$ be the iterates generated by Algorithm~\ref{alg:main.x.outline}. 
	We then define the gradient mapping and its estimator on mini-batch $\mathcal{B}$ as follows
	\begin{align}\label{def:grad_map}
	\mathcal{G}_{\eta}(x)&=\frac{1}{\eta}\left(x-\prox_{\eta\lambda\norm{\cdot}_1}(x-\eta \Grad f(x))\right),\  \text{and}\\
	\mathcal{G}_{\eta,\mathcal{B}}(x)&=\frac{1}{\eta}\left(x-\prox_{\eta\lambda\norm{\cdot}_1}(x-\eta \Grad f_\mathcal{B}(x))\right).
	\end{align}
	Here we define the noise $e(x)$ be the difference between $\mathcal{G}_{\eta}(x)$ and $\mathcal{G}_{\eta, \mathcal{B}}(x)$ with zero-mean due to the random sampling of $\mathcal{B}$, \ie, $\mathbb{E}_{\mathcal{B}}[e(x)]=0$, of which variance is bounded by $\sigma^2>0$ for one-point mini-batch.
	$\tilde{x}$ is so-called a stationary point of $F(x)$ if $\mathcal{G}_{\eta}(\tilde{x})=0$. Additionally, establishing some convergence results require the below constants to measure the least and largest magnitude of non-zero entries in $x^*$:
	\begin{align}
	0<\delta_1 :=\frac{1}{2}\min_{i\in \mathcal{I}^{\neq 0}(x^*)}|[x^*]_i|,\ \text{and}\ 0<\delta_2 :=\frac{1}{2}\max_{i\in \mathcal{I}^{\neq 0}(x^*)}|[x^*]_i|,
	\end{align}
	
	Now we state the first main theorem of \algacro{}.
	\begin{theorem}\label{thm:convex_convergence}
		Suppose $N_\mathcal{P}<\infty$ and $ N_\mathcal{O}<\infty$.
		\begin{enumerate}[label=(\roman*)]
			\item the step size $\{\alpha_k\}$ is $\mathcal{O}(1/k)$, then $\liminf_{k\to \infty}\mathbb{E}\norm{\mathcal{G}_{\alpha_k}(x_k)}_2^2=0$.
			\item $f$ is $\mu$-strongly convex, and $\alpha_k\equiv\alpha$ for any $\alpha<\min\{\frac{1}{2\mu},\frac{1}{L}\}$, then
			\begin{equation}
			\mathbb{E}[F(x_{k+1})-F^*]\leq (1-2\alpha\mu)^{\kappa_\mathcal{P}} [F(x_0)-F^*]+\frac{LC^2}{2\mu}\alpha,  
			\end{equation}
			where $\kappa_\mathcal{P}$ is the number of Prox-SG Steps employed until $k$-th iteration.
		\end{enumerate}	
	\end{theorem}
	
	Theorem~\ref{thm:convex_convergence} implies that if~\algacro{} employs~\proxsg{} Step and Orthant Step alternatively, then the gradient mapping converges to zero zero in expectation under decaying step size for general $f$ satisfying Assumption~\ref{assumption} even if $f$ is non-convex on $\mathbb{R}^n$. In other words, the iterate $\{x_k\}$ converges to some stationary point in the sense of vanishing gradient mapping. Furthermore, if $f$ is $\mu$-strongly convex and the step size $\alpha_k\equiv\alpha$ is constant, we obtain a linear convergence rate up to a solution level that is proportional to $\alpha$, which is mainly derived from the convergence properties of~\proxsg{} to optimality.  However, in practice, we may hesitate to repeatedly switch back to Prox-SG Step since most likely it is going to ruin the sparsity from the previous iterates by Orthant Step due to the stochastic nature. It is worth asking that if the convergence is still guaranteed by doing only finitely many Prox-SG Steps and then keeping doing Orthant Steps, where the below Theorem~\ref{thm:convex_convergence_plus} is drawn in line with this idea.
	
	\begin{theorem}\label{thm:convex_convergence_plus}
		Suppose $N_\mathcal{P}<\infty$, $N_\mathcal{O}=\infty$, $f$ is convex on $\{x: \norm{x-x^*}_2\leq \delta_1\}$ and $\norm{x_{N_\P}-x^*}_2\leq \frac{\delta_1}{2}$. Set $k:=N_\P+t$, $(t\in\mathbb{Z}^+)$, step size $\alpha_k=\O(\frac{1}{\sqrt{N}t})\in(0,\min\{\frac{1}{L}, \frac{\delta_1^2}{M(\delta_1+2\delta_2)}\})$, and mini-batch size $|\B_k|=\O(t)\leq N-\frac{N}{2M}$. Then for any $\tau\in (0,1)$, we have $\{x_k\}$ converges to some stationary point in expectation with probability at least $1-\tau$, \ie, $\mathbb{P}(\liminf_{k\rightarrow \infty} \mathbb{E} \norm{\mathcal{G}_{\alpha_k}(x_k)}_2^2=0)\geq 1-\tau$.
	\end{theorem}
	
	Theorem~\ref{thm:convex_convergence_plus} states the convergence is still ensured if the last iterate yielded by Prox-SG Step locates close enough to $x^*$, \ie, $\norm{x_{N_{\mathcal{P}}}-x^*}_2<\delta_1/2$. We will see in appendix that it further indicates $x^*$ inhabits the orthant faces $\{\mathcal{O}_k\}_{k\in\mathcal{S}_\mathcal{O}}$ of all subsequent iterates updated by Orthant Steps. Consequently, the convergence is then naturally followed by the property of Project Stochastic Gradient Method. Note that the local convexity-type assumption that $f$ is convex on $\{x: \norm{x-x^*}_2\leq \delta_1\}$ appears in many non-convex problem analysis, such as: tensor decomposition~\cite{ge2015escaping} and one-hidden-layer neural networks~\cite{zhong2017recovery}. Although the assumption $\norm{x_{N_{\mathcal{P}}}-x^*}_2<\delta_1/2$ is hard to be verified in practice, setting $N_\mathcal{P}$ to be large enough and $N_\mathcal{O}=\infty$ usually performs quite well, as we will see in Section~\ref{sec:exp}. To end this part, we present an upper bound of $N_\mathcal{P}$ via the probabilistic characterization to reveal that if the step size is sufficiently small, and the mini-batch size is large enough, then after $N_\mathcal{P}$~\proxsg{} Steps,~\algacro{} computes iterate $x_{N_{\mathcal{P}}}$ sufficiently close to $x^*$ with high probability. 
	
	\begin{theorem}\label{thm:n_p_upper_bound}
		Suppose $f$ is $\mu$-strongly convex on $\mathbb{R}^n$. There exists some constants $C>0, \frac{1}{2L}>\gamma>0$ such that for any constant $\tau\in(0,1)$, if $\alpha_k$ satisfies $\alpha_k\equiv\alpha<\min\left\{\frac{2\gamma\mu\tau\delta_1^2}{(2L\gamma-1){C}}, \frac{1}{2\mu},\frac{1}{L}\right\}$, and the mini-batch size $|\mathcal{B}_k|$ satisfies $|\mathcal{B}_k|> \frac{8\gamma\mu \sigma^2}{2\gamma\mu\tau\delta_1^2-(2L\gamma-1)C\alpha}$, then the probability of $\norm{x_{N_\mathcal{P}}-x^*}_2\leq \delta_1/2$ is at least $1-\tau$ for any $N_{\mathcal{P}}\geq K$ where $K:=\left\lceil\frac{\log{(\text{poly}(\tau\delta_1^2, 1/|\mathcal{B}_k|, \alpha)/(F(x_0)-F^*))}}{\log{(1-2\mu\alpha)}}\right\rceil$ and $\text{poly}(\cdot)$ represents some polynomial of $\tau\delta_1^2, 1/|\mathcal{B}_k|$ and $\alpha$.
	\end{theorem}
	In words, Theorem~\ref{thm:n_p_upper_bound} implies that after sufficient number of iterations, with high probability Prox-SG produces an iterate $x_{N_\mathcal{P}}$ that is $\delta$-close to $x^*$. However, we note that it does not guarantee $x_{N_\mathcal{P}}$ as sparse as $x^*$;  as we explained before, due to the limited projection region and randomness, $x_{N_\mathcal{P}}$ may still have a large number of non-zero elements, though many of them could be small. As will be demonstrated in Section~\ref{sec:exp}, the following Orthant Steps will significantly promote the sparsity of the solution.

	\section{Numerical Experiments}\label{sec:exp}
	
	In this section, we consider solving $\ell_1$-regularized classification tasks with both convex and non-convex approaches. In Section~\ref{sec:convex_exp}, we focus on logistic regression (convex), and compare~\algacro{} with other state-of-the-art methods including~\proxsg,~\rda{} and~\proxsvrg{} on numerous datasets. Three evaluation metrics are used for comparison: (i) final objective function value, (ii) density of the solution (percentage of nonzero entries), and (iii) runtime. Next, in Section~\ref{sec:non-convex_exp}, we apply~\algacro{} to deep neural network (non-convex) with popular architectures designed for classification tasks to further demonstrate its effectiveness and superiority. For these extended non-convex experiments, we also evaluate the generalization performance on unseen test data.

	\subsection{Convex setting: logistic regression }\label{sec:convex_exp}
	We first focus on the convex $\ell_1$-regularized logistic regression with the form
	\begin{equation}\label{def:minimize_logistic_l1}
	\minimize{(x; b)\in \R^{n+1}}\ \frac{1}{N}\sum_{i=1}^N \log(1 + e^{-l_i (x^T d_i +b)}) + \lambda \norm{x}_1,
	\end{equation}
	for binary classification, where $N$ is the number of samples, $n$ is the feature size of each sample, $b$ is the bias, $d_i \in \R^n$ is the vector representation of the $i$-th sample, $l_i \in \{-1, 1\}$ is the label of the $i$-th sample, and $\lambda$ is the regularization parameter. We set  $\lambda=1/N$ throughout the convex experiments, and  test~problem~\eqref{def:minimize_logistic_l1} on 8 public large-scale datasets from LIBSVM repository~\footnote{\url{https://www.csie.ntu.edu.tw/~cjlin/libsvmtools/datasets/}}, as summarized in Table~\ref{table:datasets}. 
	
	\begin{table}[t]
		\centering
		\def\arraystretch{1.1}
		\caption{Summary of datasets\label{table:datasets}}
		\resizebox{\textwidth}{!}{
			\begin{tabular}{@{\extracolsep{4pt}}ccccc?cccc}
				\Xhline{2\arrayrulewidth}
				Dataset & N & n  & Attribute & & Dataset & N & n  & Attribute \\
				\hline
				a9a & 32561 & 123 & binary \{0, 1\} & & real-sim & 72309 & 20958 & real [0, 1]\\
				higgs & 11000000 & 28 & real $[-3, 41]$ & & rcv1 & 20242 & 47236 &  real [0, 1]\\
				kdda & 8407752 & 20216830 & real $[-1, 4]$ &  & url\_combined & 2396130 & 3231961 & real $[-4, 9]$\\
				news20 & 19996 & 1355191 &   unit-length&  & 	w8a & 49749 & 300   & binary \{0, 1\}\\
				\Xhline{2\arrayrulewidth}
			\end{tabular}
		}
		
		
	\end{table}
	
	We train the models with a maximum number of epochs as $30$. Here “one epoch” means we partition $\{1,\cdots,N\}$ uniformly at random into a set of mini-batches. The mini-batch size $|\mathcal{B}|$ for all the convex experiments is set to be $\min\{256, \ceil*{0.01N}\}$ similarly to~\cite{yang2019stochastic}. The step size $\alpha_k$ for \proxsg{}, \proxsvrg{} and \algacro{} is initially set to be $1.0$, and decays  every epoch with a factor $0.995$. For \rda,  we fine tune its hyperparameter $\gamma$ per dataset to reach the best results. The switching between Prox-SG Step and~\orthantstep{} plays a crucial role in OBProx-SG. Following Theorem~\ref{thm:convex_convergence}(i), we set $N_\mathcal{P}=N_\mathcal{O}=5N/|\mathcal{B}|$ in Algorithm~\ref{alg:main.x.outline}, namely first train the models $5$ epochs by~\objstep, followed by performing~\orthantstep{} $5$ epochs, and repeat such routine until the maximum number of epochs is reached. Inspired by Theorem~\ref{thm:convex_convergence}(ii), we also test~\algacroplus{} with $N_\mathcal{P}=15N/|\mathcal{B}|, N_\mathcal{O}= \infty$ such that after 15 epochs of Prox-SG Steps we stick to Orthant Step till the end. Experiments are conducted on a 64-bit machine with an 3.70GHz Intel Core i7 CPU and 32 GB of main memory.
	
\begin{figure}[t]
		\centering
		\includegraphics[width=\textwidth]{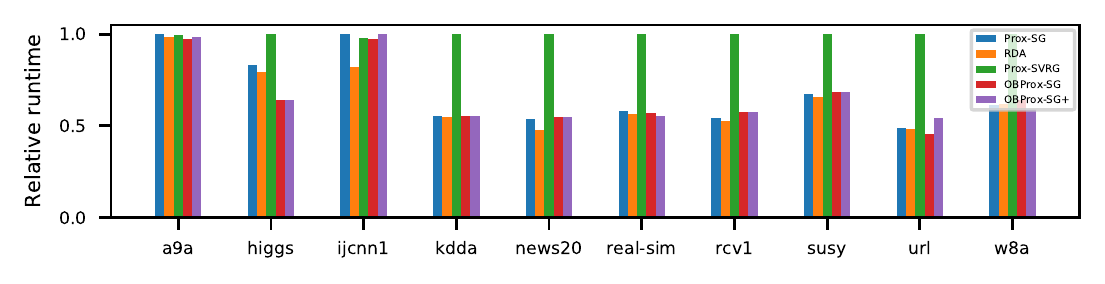}
		\caption{Relative runtime for tested algorithms on convex problems}
		\label{figure:runtime_convex}
	\end{figure}

	We compare the performance of~\algacroall{} with other methods on the datasets in Table~\ref{table:datasets}, and report the final objective value $F$ and $f$ (Table~\ref{table:f_value_convex}), density (percentage of non-zero entries) in the solution (Table~\ref{table:density_convex}) and runtime (Figure~\ref{figure:runtime_convex}). For ease of comparison, we mark the best result as bold in the tables. 
	
	Our observations are summarized as follows. Table~\ref{table:f_value_convex} shows that our~\algacroall{} performs significantly better than~\rda{}, and is competitive to~\proxsg{} and~\proxsvrg{} in terms of the final $F$ and $f$ (round up to 3 decimals), which implies that~\algacroall{},~\proxsg{} and~\proxsvrg{} can reach comparable convergence results in practice. Besides the convergence, we have a special concern about the sparsity of the solutions. As is demonstrated in Table~\ref{table:density_convex}, \algacroall{}
	is no doubt the best solver. In fact, \algacro{} achieves the solutions of highest sparsity (lowest density) on 1 out of 8 datasets, while \algacroplus{} performs even better, which computes all solutions with the highest sparsity. Apparently,  \algacroall{} has strong superiority in promoting sparse solutions while retains almost the same accuracy. Finally, for runtime comparison, we plot the relative runtime of these solvers, which is scaled by the maximum runtime consumed by a particular solver on that dataset. Figure~\ref{figure:runtime_convex} indicates that \proxsg,~\rda{} and \algacroall{} are almost as efficient as each other, while~\proxsvrg{} takes much more time due to the computation of full gradient.
	
	\begin{table}[t]
		\centering
		\caption{Objective function values $F/f$ for tested algorithms on convex problems}
		\label{table:f_value_convex}
		\def\arraystretch{1.1}
		
		{\scriptsize
			\begin{tabularx}{\textwidth} { cccccc
				}
				\Xhline{3\arrayrulewidth}
				
				Dataset & \proxsg{} & \rda{} & \proxsvrg{} & \algacro{} & \algacroplus{}\\
				\hline
				a9a & \ 0.332\ /\ 0.330 & \ 0.330\ /\ 0.329 & \ 0.330\ /\ 0.329 & \ \textbf{ 0.327}\ /\ \textbf{0.326} & \ {0.329}\ /\ { 0.328}  \\ 
				higgs & \ \textbf{0.326}\ /\ \textbf{0.326} & \ \textbf{0.326}\ /\ \textbf{0.326} & \ \textbf{0.326}\ /\ \textbf{0.326} & \ \textbf{0.326}\ /\ \textbf{0.326} & \ \textbf{0.326}\ /\ \textbf{0.326} \\
				kdda & \ \textbf{0.102}\ /\ \textbf{0.102} & \ 0.103\ /\ 0.103 & \ 0.105 \ /\ 0.105  & \ \textbf{0.102}\ /\ \textbf{0.102} & \ \textbf{0.102}\ /\ \textbf{0.102}\\
				news20 & \ \textbf{0.413}\ /\ \textbf{0.355} & \ 0.625\ /\ 0.617 & \ \textbf{0.413}\ /\ \textbf{0.355} & \ \textbf{0.413}\ /\ \textbf{0.355} &\  \textbf{0.413}\ /\ \textbf{0.355}  \\ 
				real-sim & \ \textbf{0.164}\ /\ \textbf{0.125} & \ 0.428\ /\ 0.421 & \ \textbf{0.164}\ /\ \textbf{0.125} & \ \textbf{0.164}\ /\ \textbf{0.125} & \ \textbf{0.164}\ /\ \textbf{0.125}  \\ 
				rcv1 & \ \textbf{0.242}\ /\ \textbf{0.179} & \ 0.521\ /\ 0.508 & \ \textbf{0.242}\ /\ \textbf{0.179}  & \ \textbf{0.242}\ /\ \textbf{0.179}  & \ \textbf{0.242}\ /\ \textbf{0.179}  \\ 
				url\_combined & \ {0.050\ /\  0.049} & \ 0.634\ /\ 0.634  & \ 0.078\ /\ 0.077 & \ { 0.050\ /\  0.049} & \ \textbf{0.047}\ /\ \textbf{0.046}\\
				w8a & \ \textbf{0.052}\ /\	\textbf{0.048} & \ 0.080\ /\	0.079 & \ \textbf{0.052}\ /\	\textbf{0.048} & \ \textbf{0.052}\ /\	\textbf{0.048} & \ \textbf{0.052}\ /\	\textbf{0.048}  \\
				\Xhline{3\arrayrulewidth}
			\end{tabularx}
		}
		\caption{Density (\%) of solutions for tested algorithms on convex problems}
		\label{table:density_convex}
		{\scriptsize
			\begin{tabularx}{\textwidth} { 
					>{\centering\arraybackslash}X 
					>{\centering\arraybackslash}X 
					>{\centering\arraybackslash}X 
					>{\centering\arraybackslash}X 
					>{\centering\arraybackslash}X 
					>{\centering\arraybackslash}X }
				\Xhline{3\arrayrulewidth}
				Dataset & \proxsg{} & \rda{} & \proxsvrg{} & \algacro{} & \algacroplus{} \\
				\hline
				a9a & 96.37 & 86.69 & {61.29} & 62.10 & \textbf{59.68}  \\ 
				higgs & 89.66 & 96.55 & 93.10 & \textbf{70.69} & \textbf{70.69} \\
				kdda & 0.09 & 18.62 & 3.35 & {0.08} & \textbf{0.06}\\
				news20 & 4.24 & 0.44 & {0.20} & {0.20} & \textbf{0.19}   \\ 
				real-sim & 53.93 & 52.71 & {22.44} & {22.44} & \textbf{22.15}   \\ 
				rcv1 & 16.95 & 9.61 & {4.36} & {4.36} & \textbf{4.33}   \\ 
				url\_combined & 7.73 & 41.71 & 6.06 & {3.26} & \textbf{3.00} \\
				w8a & 99.00 & 99.83 & 78.07 & {78.03} & \textbf{74.75}  \\ 
				\Xhline{3\arrayrulewidth}
			\end{tabularx}
		}
	\end{table}

	The above experiments in convex setting demonstrate that the proposed \algacroall{} outperform the other state-of-the-art methods, and have apparent strengths in generating much sparser solutions efficiently and reliably. 
	
	\subsection{Non-convex setting: deep neural network}\label{sec:non-convex_exp}
	
	We now apply~\algacroall{} to the non-convex setting that solves classification tasks by Deep Convolutional Neural Network (CNN) on the benchmark datasets CIFAR10~\cite{Krizhevsky09} and~\fashionmnist{}~\cite{xiao2017online}. Specifically, we are testing two popular CNN architectures, \ie,  \mobilenet{}~\cite{howard2017mobilenets} and \resnet{}~\cite{he2016deep}, both of which have proven successful in many image classification applications. We add an $\ell_1$-regularization term to the raw problem, where $\lambda$ is set to be $10^{-4}$ throughout the non-convex experiments.
	
	
	
	We conduct all non-convex experiments for 200 epochs with a mini-batch size of 128 on one GeForce GTX 1080 Ti GPU. The step size $\alpha_k$ in \proxsg{}, \proxsvrg{} and \algacroall{} is initialized as $0.1$, and decay by a factor 0.1 periodically. The $\gamma$ in \rda{} is fine-tuned to be $20$ for \cifar{} and $30$ for \fashionmnist{} in order to achieve the best performance. Similar to convex experiments, we set $N_\mathcal{P}=N_\mathcal{O}=5N/|\mathcal{B}|$ in~\algacro, and set $N_\mathcal{P}=100N/|\mathcal{B}|$, $N_\mathcal{O}=\infty$ in~\algacroplus{} since running~\proxsg{} Step 100 epochs already achieves an acceptable validation accuracy.  
	
	Based on the experimental results, the conclusions that we  made  previously in convex setting still hold in the current non-convex case: (i) \algacroall{} performs competitively among the methods with respect to the final objective function values, see Table~\ref{table:object_value_nonconvex}; (ii) \algacroall{} computes much sparser solutions which are significantly better than other methods as shown in Table~\ref{table:density_nonconvex}. Particularly,~\algacroplus{} achieves the highest sparse (lowest dense) solutions on all non-convex tests, of which the solutions are 4.24 to 21.86 times sparser than those of~\proxsg{}, while note that \rda{} and~\proxsvrg{} perform not comparable on the sparsity exploration because of the ineffectiveness of variance reduction techniques for deep learning~\cite{defazio2019ineffectiveness}. In addition, we evaluate how well the solutions generalize on unseen test data. Table~\ref{table:density_nonconvex} shows that all the methods reach a comparable testing accuracy except RDA.

	\begin{table}[t]
		\centering
		\caption{Final objective values $F/f$  for tested algorithms on non-convex problems}
		\label{table:object_value_nonconvex}
		\def\arraystretch{1.1}
		\resizebox{\textwidth}{!}{
			\begin{tabular}{@{\extracolsep{4pt}}ccccccc @{}}
				\Xhline{3\arrayrulewidth}
				Backbone & Dataset &\proxsg{} & \rda{} & \proxsvrg{} & \algacro{} & \algacroplus{}\\
				\hline
				\multirow{2}{*}{\mobilenet{}}& \cifar{} & {1.473}\ /\ {0.049} & 4.129\ /\ 0.302 & 1.921\ /\ 0.079 & 1.619\ /\	\textbf{0.048} & \textbf{1.453}\ /\ 0.063\\ 
				& \fashionmnist{} & {1.314}\ /\ \textbf{0.089} & 4.901\ /\ 0.197 & 1.645\ /\ 0.103 & 2.119\ /\ \textbf{0.089} & \textbf{1.310}\ /\	{0.099}\\ \hdashline
				\multirow{2}{*}{\resnet{}} & \cifar{} & 0.781\ /\ 0.034 & 1.494\ /\	0.051 &	0.815\ /\ {0.031} & 	\textbf{0.746}\ /\ \textbf{0.021} & {0.755}\ /\ 0.044\\
				&\fashionmnist{} & 0.688\ /\ 0.103& 1.886\ /\	{0.081} & {0.683}\ /\ \textbf{0.074} & \textbf{0.682}\ /\ \textbf{0.074} & 0.689\ /\ 0.116\\
				\Xhline{3\arrayrulewidth} 
			\end{tabular}
		}
		{\\}
		\caption{Density/testing accuracy (\%) for tested algorithms on non-convex problems}
		\label{table:density_nonconvex}
		{\scriptsize
			\begin{tabularx}{\textwidth} {  ccccccc
				}
				\Xhline{3\arrayrulewidth}
				Backbone & Dataset &\proxsg{} & \rda{} & \proxsvrg{} & \algacro{} & \algacroplus{}\\
				\hline
				\multirow{2}{*}{\mobilenet{}}& \cifar{} & 14.17/\textbf{90.98} & 74.05/81.48 & 92.26/87.85 & 9.15/90.54 & \textbf{\textbf{2.90}}/90.91\\
				& \fashionmnist{} & 5.28/94.23 & 74.67/92.12 & 75.40/93.66 & 4.15/94.28 & \textbf{\textbf{1.23}}/\textbf{94.39}\\\hdashline
				\multirow{2}{*}{\resnet{}} & \cifar{} & 11.60/92.43 & 41.01/90.74 & 37.92/92.48 & 2.12/\textbf{92.81} & \textbf{\textbf{0.88}}/92.45 \\
				&\fashionmnist{} & 6.34/94.28 & 42.46/93.66 & 35.07/94.24 & 5.44/\textbf{94.39} & \textbf{\textbf{0.29}}/93.97\\
				\Xhline{3\arrayrulewidth} 
			\end{tabularx}
		}
		
	\end{table}

	\begin{figure}[t]
		\includegraphics[width=0.48\textwidth]{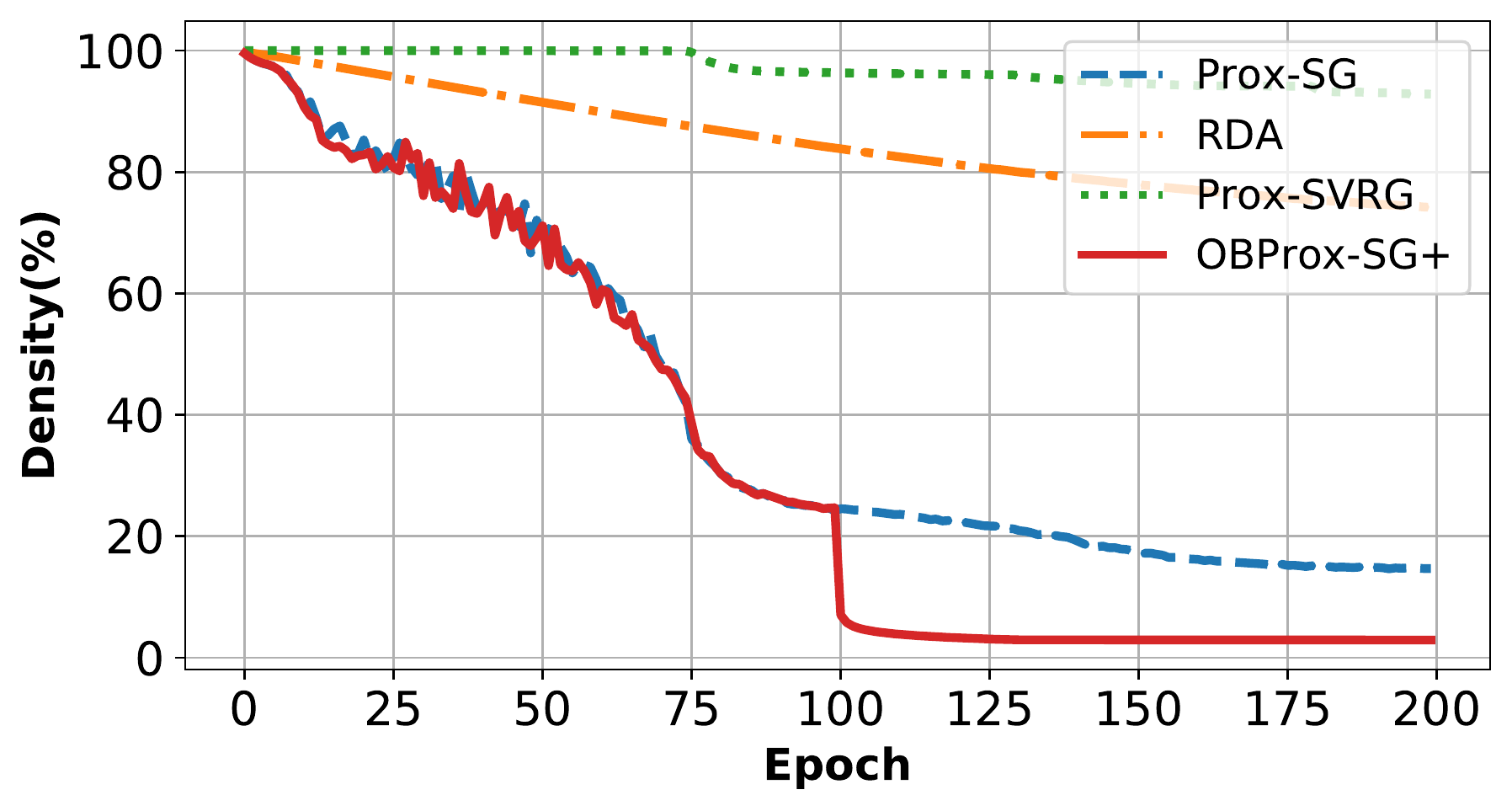}
		\includegraphics[width=0.48\textwidth]{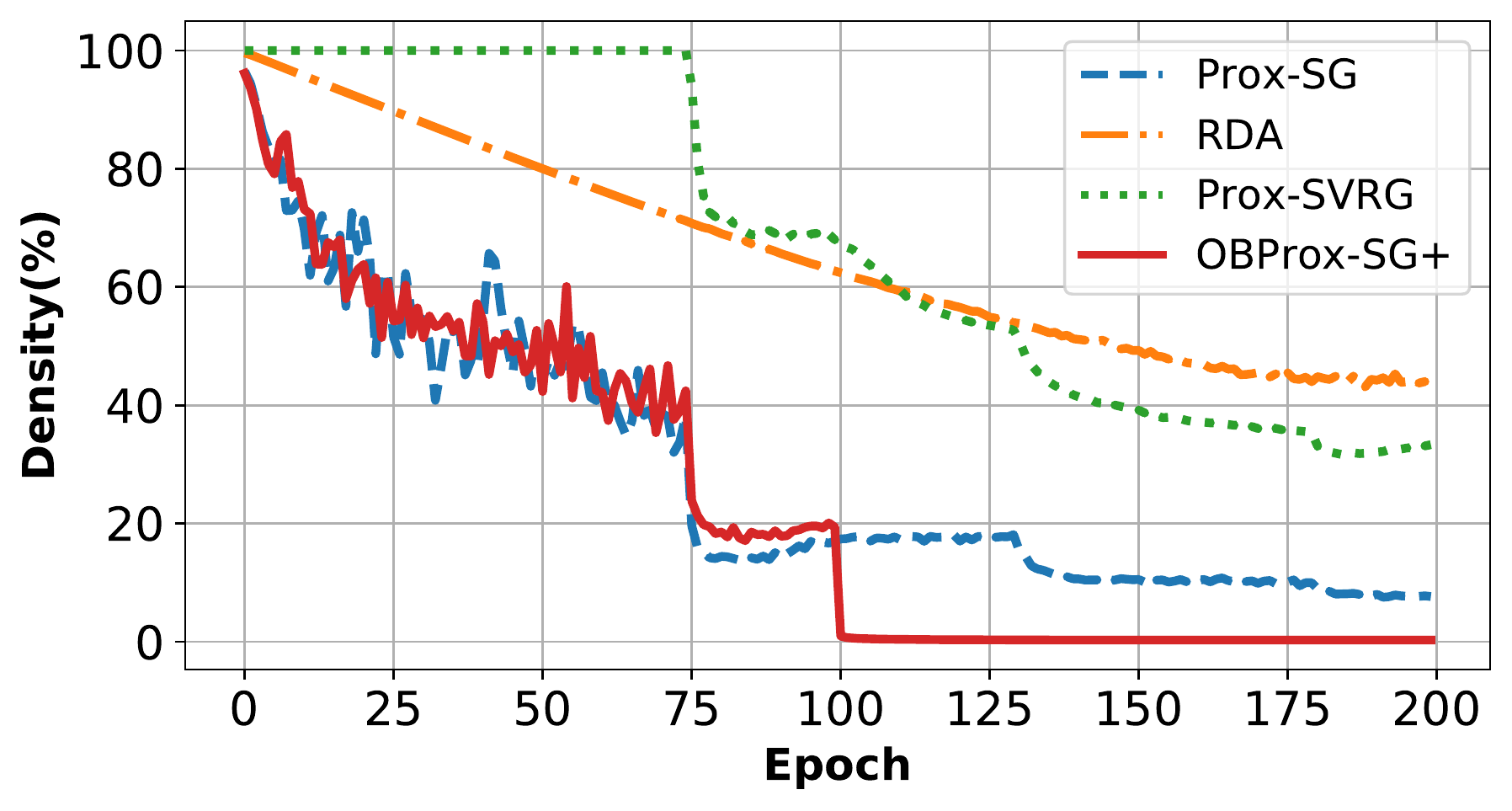}
		\caption{ Density. (L): \mobilenet{} on \cifar{}. (R): \resnet{} on \fashionmnist{}}
		\label{figure:density_non-convex}
	\end{figure}

	Finally, we investigate the sparsity evolution of the iterates to reveal the superiority of~\orthantstep{} on sparsity promotion, where we use~\algacroplus{} as the representative of~\algacroall{} for illustration. As shown in Figure~\ref{figure:density_non-convex}, \algacroplus{} produces the 
	highest sparse (lowest dense) solutions compared with other methods. Particularly, at the early $N_{\mathcal{P}}$ iterations, ~\algacroplus{} performs merely the same as~\proxsg{}. However, after the switching to \orthantstep{} at the 100th epoch, \algacroplus{} outperforms all the other methods dramatically. It is a strong evidence that because of  the construction of orthant face subproblem and the larger projection region, our orthant based technique is more remarkable than the standard proximal gradient step and its variants in terms of the sparsity exploration. As a result, the solutions computed by~\algacro{} generally have a better interpretation under similar generalization performances. Furthermore,~\algacro{} may be further used to save memory and hard disk storage consumption drastically by constructing sparse network architectures.

	\section{Conclusions}
	
	We proposed an Orthant Based Proximal Stochastic Gradient Method (OBProx-SG) for solving $\ell_1$-regularized problem, which combines the advantages of deterministic orthant based methods and proximal stochastic gradient method. 
	In theory, we proved that it converges to some global solution in expectation for convex problems and some stationary point for non-convex formulations. Experiments on both convex and non-convex problems demonstrated that OBProx-SG usually achieves competitive objective values and much sparser solutions compared with state-of-the-arts stochastic solvers. 
	
	\section*{Acknowledgments}
	We would like to thank 
	the four anonymous reviewers for their constructive comments.
	T. Ding was partially supported by NSF grant 1704458. Z. Zhu was partially supported by NSF grant 2008460.

	%
	%
	%
	\bibliographystyle{splncs04}

\end{document}